\documentclass[12pt,oneside]{article}

\usepackage[english]{babel}
\usepackage[T1]{fontenc}
\usepackage{graphicx}
\usepackage{amsfonts}
\usepackage{amsmath}
\usepackage{amssymb}
\usepackage{mathrsfs}
\usepackage{amsthm}
\usepackage{latexsym}
\usepackage[all]{xy}
\usepackage{tikz}
\usepackage{tkz-arith}
\usetikzlibrary{arrows,shapes,matrix,decorations.pathmorphing,shapes.geometric,petri}

\usepackage{float}                               
\restylefloat{figure}                              

\newtheorem{teo}{Theorem}[]

\theoremstyle{definition}

\DeclareMathOperator{\rc}{\stackrel{\rightarrow}{rc}}

\begin{document}

\centerline{\bf A Note on the Rainbow Connectivity of Tournaments}
\vskip 1pc
\centerline{Jes\'us Alva-Samos\footnote{
              Instituto de Matem\'aticas, UNAM}  \hskip 1pc
        Juan Jos\'e Montellano-Ballesteros\footnote{Instituto de Matem\'aticas, UNAM. 
                juancho@matem.unam.mx }}
\vskip 1pc             
\centerline{April 27, 2015}

\begin{abstract}
An arc-coloured digraph $D$ is said to be \emph{rainbow connected} if for every two vertices $u$ and $v$ there is an $uv$-path all whose arcs have different colours. The minimun number of colours required to make the digraph rainbow connected is called the \emph{rainbow connection number} of $D$, denoted $\rc(D)$.  In \cite{Dorbec} it was showed that if $T$ is a strong tournament with $n\ge5$ vertices, then $2\le\rc(T)\le n-1$; and that 
for every $n$ and $k$ such that $3\le k\le n-1$, there exists a tournament $T$ on $n$ vertices such that $\rc(T)=k$.

In this note it is showed that for any $n\ge6$, there is a tournament $T$ of $n$ vertices such that $\stackrel{\rightarrow}{rc}(T)=2$.
\end{abstract}


\section{Introduction}
Chartrand, Johns, McKeon and Zhang introduced in \cite{Ch} the concept of rainbow connection in graphs, and since then such topic has been broadly studied, see \cite{Li} for a survey on the matter. The rainbow connection for digraphs was first presented by Dorbec, Schiermeyer, Sidorowicz and Sopena in \cite{Dorbec}.

A \emph{tournament} is a digraph where every two vertices has exactly one arc joining them. Let $\rho$ be an arc-colouring on a digraph $D$, a (directed) path in $D$ is called a \emph{rainbow} if it has no two arcs sharing the same colour. If there is a rainbow between every pair of vertices, then the digraph is \emph{rainbow connected} and $\rho$ is called a \emph{rainbow colouring}. The \emph{rainbow connection number} of a digraph $D$, denoted by $\rc(D)$, is defined as the smallest number of colours such that $D$ admits a rainbow colouring.

In \cite{Dorbec} the following two theorems were proven:

\begin{teo}[Dorbec et al. \cite{Dorbec}]
If $T$ is a strong tournament with $n\ge5$ vertices, then $2\le\rc(T)\le n-1$.
\end{teo}

\begin{teo}[Dorbec et al. \cite{Dorbec}]\label{Texists}
For every $n$ and $k$ such that $3\le k\le n-1$, there exists a tournament $T$ on $n$ vertices such that $\rc(T)=k$.
\end{teo}

Also, in \cite{Dorbec}, it is showed that for each $n\cong 8$ mod $12$ there is a tournament  $T$ of order $n$ and $\rc(T)=2$. 

In this note it is showed that for every $n\geq 6$ there is a tournament  $T$ of order $n$ and $\rc(T)=2$.

\section{The tournaments}
When $n=4,5$ it is easy to verified that $\rc(T)\ge3$ for each tournament $T$ on $n$ vertices.   For an integer $n\geq2$ and a set $S\subseteq\{1,2,\dots,n-1\}$, the \emph{circulant digraph} $C_n(S)$ is defined as the digraph with vertex set $V(C_n(S))=\{v_0,v_1,\dots,v_{n-1}\}$ and arc set $A(C_n(S))=\{(v_i,v_j)\mid j-i\stackrel{n}{\equiv}s,\ s\in S\}$.

\begin{teo}\label{Texistsb}
For every $n\ge6$, there is a tournament $T$ of order $n$ with $\rc(T)=2$.
\end{teo}

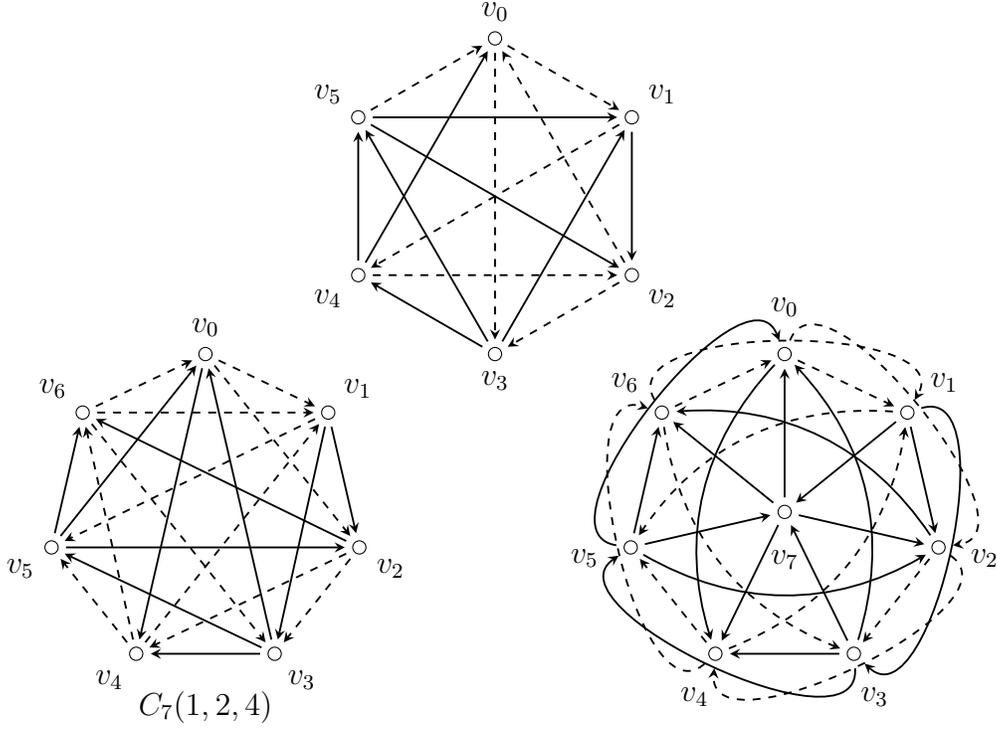
\begin{figure}[ht]
\begin{center}
\begin{tikzpicture}[every circle node/.style ={circle,draw,minimum size= 5pt,inner sep=0pt, outer sep=0pt},
every rectangle node/.style ={}];

\begin{scope}[scale=0.7, xshift=5.5cm, yshift=6cm]
\node [circle] (0) at (90:3)[label=90:$v_0$]{};	\node [circle] (1) at (30:3)[label=30:$v_1$]{};
\node [circle] (2) at (330:3)[label=330:$v_2$]{};	\node [circle] (3) at (270:3)[label=270:$v_3$]{};
\node [circle] (4) at (210:3)[label=210:$v_4$]{};	\node [circle] (5) at (150:3)[label=150:$v_5$]{};
\foreach \from/\to in {0/1,0/3,1/4,2/0,2/3,4/2,5/0}
\draw [->, shorten <=3pt, shorten >=3pt, >=stealth, line width=.7pt, dashed] (\from) to  (\to);
\foreach \from/\to in {1/2,3/1,3/4,3/5,4/5,4/0,5/1,5/2}
\draw [->, shorten <=3pt, shorten >=3pt, >=stealth, line width=.7pt] (\from) to  (\to);
\end{scope}

\begin{scope}[scale=0.7, xshift=0cm, yshift=0cm]
\node [circle] (0) at (90:3)[label=90:$v_0$]{};	\node [circle] (1) at (39:3)[label=40:$v_1$]{};
\node [circle] (2) at (347:3)[label=350:$v_2$]{};	\node [circle] (3) at (296:3)[label=290:$v_3$]{};
\node [circle] (4) at (244:3)[label=250:$v_4$]{};	\node [circle] (5) at (193:3)[label=190:$v_5$]{};
\node [circle] (6) at (141:3)[label=140:$v_6$]{};	\node (99) at (270:3.7){$C_7(1,2,4)$};
\foreach \from/\to in {0/1,0/2,1/5,2/3,2/4,4/5,4/6,4/1,6/0,6/1,6/3}
\draw [->, shorten <=3pt, shorten >=3pt, >=stealth, line width=.7pt,dashed] (\from) to  (\to);
\foreach \from/\to in {0/4,1/2,1/3,2/6,3/4,3/5,3/0,5/6,5/0,5/2}
\draw [->, shorten <=3pt, shorten >=3pt, >=stealth, line width=.7pt] (\from) to  (\to);
\end{scope}

\begin{scope}[scale=0.7, xshift=11cm, yshift=0cm]
\node [circle] (0) at (90:3)[]{};	\node [rectangle] (00) at (90:3.9)[]{$v_0$};
\node [circle] (1) at (39:3)[]{};	\node [rectangle] (00) at (39:3.9)[]{$v_1$};
\node [circle] (2) at (347:3)[]{};	\node [rectangle] (00) at (347:3.9)[]{$v_2$};
\node [circle] (3) at (296:3)[]{};	\node [rectangle] (00) at (296:3.9)[]{$v_3$};
\node [circle] (4) at (244:3)[]{};	\node [rectangle] (00) at (244:3.9)[]{$v_4$};
\node [circle] (5) at (193:3)[]{};	\node [rectangle] (00) at (193:3.9)[]{$v_5$};
\node [circle] (6) at (141:3)[]{};	\node [rectangle] (00) at (141:3.9)[]{$v_6$};
\node [circle] (7) at (0,0) []{};	\node [rectangle] (00) at (270:0.9)[]{$v_7$};
\foreach \from/\to in {0/1,2/3,4/5,6/0}
\draw [->, shorten <=3pt, shorten >=3pt, >=stealth, line width=.7pt,dashed] (\from) to  (\to);
\foreach \from/\to in {1/2,3/4,5/6,  7/0,7/2,7/4,7/6,1/7,3/7,5/7}
\draw [->, shorten <=3pt, shorten >=3pt, >=stealth, line width=.7pt] (\from) to  (\to);
\foreach \from/\to in {0/2,2/4,4/6,6/1}
\draw [->,shorten <=3pt,shorten >=3pt,>=stealth,line width=.7pt,dashed] (\from) to [bend left=120] (\to);
\foreach \from/\to in {1/3,3/5,5/0}
\draw [->,shorten <=3pt,shorten >=3pt,>=stealth,line width=.7pt] (\from) to [bend left=120] (\to);
\foreach \from/\to in {1/5,4/1,6/3}
\draw [->,shorten <=3pt,shorten >=3pt,>=stealth,line width=.7pt,dashed] (\from) to [bend right=30] (\to);
\foreach \from/\to in {0/4,2/6,3/0,5/2}
\draw [->,shorten <=3pt,shorten >=3pt,>=stealth,line width=.7pt] (\from) to [bend right=30] (\to);
\end{scope}

\end{tikzpicture}
\caption{Tournaments on $n$ vertices with $\rc(T_n)=2$ where $n=6,7,8$.}\label{fig:trnmnts}
\end{center}
\end{figure}

\begin{proof}
For $n=6$, let $T$ be the first tournament in Figure \ref{fig:trnmnts}.
Now, for $n=2k+1$, with $k\ge3$, let $T=C_{2k+1}(1,2,4,\dots,2(k-1))$ and consider the partition of $A(T)$ into the sets
\begin{eqnarray*}
A_0 & = & \{(u_0,u_1),(u_0,u_2),(u_1,u_{2k-1})\} \\ 
   &  & \cup\big(\{(u_r,u_s)\mid r\equiv0\text{ mod }2, r\ge2\}\setminus\{(u_2,u_{2k})\}\big)
\end{eqnarray*}
and $A_1=A(T)\setminus A_0$. Now, let $\rho$ be the coloring where $\rho(a)=i$ if $a\in A_i$ for $i=0,1$.  We will show that  there is a rainbow $u_iu_j$-path for $0\le i\le n-1$ and $j=i+3,i+5,\dots,i+2k-1,i+2k$, and the result will follow.

For each $i$ with $3\leq i \leq 2k-1$ let consider the paths 
$$\left\{\begin{array}{ll} u_iu_{i+1}u_r & \text{for  $r=3, 5, \dots,2k-1$}; \\ u_iu_{i+2k-2}u_{i+2k} ; \end{array}\right.$$
and for the vertex $u_{2k}$ consider the paths $u_{2k} u_1 u_2$,  $u_{2k} u_{2k-3} u_{2k-1}$ and  for $r=5, 7,\dots 2k-1$ let $u_{2k} u_0 u_r$. 

For $i\in\{1,2\}$  consider the paths $u_iu_{i+1}u_{i+r}$ for $r= 3, 5,\dots 2k-3$; and the paths  $u_1u_{2k-1}u_{2k}$,   $u_1u_{2k-1}u_{0}$,  $u_2u_{2k}u_{0}$, $u_2u_{2k}u_{1}$.

Finally,  for the vertex $u_0$ let consider the paths 
$$\left\{\begin{array}{ll} u_0u_1u_3 &\\
u_0u_{r-1}u_r & \text{for  $r=5, 7,\dots 2k-1$ }; \\ u_0u_{2k-2}u_{2k}.\end{array}\right.$$

By the definition of $\rho$, it is not hard to see that all these paths are rainbow paths and therefore $\rc(T)=2$.

  For the case $n=2k\ge8$, let consider  the tournament $T$ obtained by adding a new vertex $v$ to $C_{2k-1}(1,2,4,\dots,2(k-2))$ such that $(v,u_i)\in A(T)$ if  and only if   $i$ is even.  Let colouring the arcs of $T$ by extending the coloring $\rho$ of $C_{2k-1}(1,2,4,\dots,2(k-2))$ described above and given color 1 to  all the new arcs.   Consider the paths $vu_iu_{i+1}$ for $i+1$ odd;  $u_iu_{i+1}v$ for $i$ even, $i\not=2k$; and $u_{2k}u_1v$.  
 By construction all these paths are rainbow paths, therefore $\rc(T)=2$.
\end{proof}

Finally, combining Theorems \ref{Texists} and \ref{Texistsb} we have

\begin{teo}
For every $n\ge6$ and every $k$ such that $2\le k\le n-1$, there exists a tournament $T$ on $n$ vertices such that $\rc(T)=k$.
\hfill{$\square$}
\end{teo}

\end{document}